\documentclass[11pt]{article} 
\usepackage{amsmath,amsthm,amscd,amssymb}
\usepackage{latexsym}
\usepackage{epsf}
\usepackage{epsfig}

\newcommand{\diff}{\operatorname{Diff^2}}
\newcommand{\Diff}{\operatorname{Diff^2}}

\newcommand{\vol}{\operatorname{Vol}}
\renewcommand{\div}{\operatorname{div}}
\newcommand{\R}{\mathbb{R}}

\newcommand{\w}{\omega}
\newcommand{\om}{\omega}

\newcommand{\ep}{\varepsilon}
\newcommand{\al}{\alpha}
\newcommand{\la}{\lambda}

\newcommand{\SM}{\mathcal{M}}
\newcommand{\SP}{\mathcal{P}}
\newcommand{\Chi}{\mathfrak{X}}
\newcommand{\J}{distortion minimal }

\theoremstyle{plain}
\newtheorem{theorem}{Theorem}[section]
\newtheorem{lemma}[theorem]{Lemma}
\newtheorem{corollary}[theorem]{Corollary}
\newtheorem{proposition}[theorem]{Proposition}

\theoremstyle{definition}
\newtheorem{definition}[theorem]{Definition}

\newtheorem{example}[theorem]{Example}
\theoremstyle{remark}
\newtheorem{remark}[theorem]{Remark}

\newtheorem{case[theorem]}{Case}

\title{Distortion Minimal Morphing I: \\The Theory For Stretching}
\author{Oksana Bihun\footnote{E-mail: oksana@math.missouri.edu}\hspace*{.05in} and Carmen Chicone\footnote{Corresponding Author E-mail: carmen@math.missouri.edu}\\Department of Mathematics\\University of
Missouri-Columbia\\Columbia, Missouri 65211, USA}

\begin{document}
\maketitle
\begin{abstract}
We consider the problem of distortion minimal morphing of $n$-dimensional 
compact connected oriented smooth manifolds without boundary embedded in 
$\R^{n+1}$. Distortion involves bending and stretching. In this paper,  
minimal distortion (with respect to stretching) is defined as the 
infinitesimal relative change in volume. The existence of minimal 
distortion diffeomorphisms between diffeomorphic manifolds is proved.  
A definition of minimal distortion morphing between two isotopic manifolds 
is given, and  the existence of 
minimal distortion morphs between every pair of isotopic embedded manifolds 
is proved.
\end{abstract}

optimization

\section{Introduction}
A morph is a transformation between two shapes through a set of
intermediate shapes.  A minimal morph is such a transformation that minimizes distortion.  

There are important applications of minimal morphing in  
manufacturing~\cite{LYS, YPM},
computer graphics~\cite{W,W1}, movie making~\cite{HLZ}, 
and mesh construction~\cite{HH,LDSS}. To address these applications would require a theory of minimal morphing that includes bending and stretching together with algorithms to compute minimal morphs. In this paper, we do not address applications; rather, we formulate and solve the mathematical problem of minimal morphing with respect to stretching. A complete theory of minimal morphing is a subject for future research.      

In the research literature available at present,  minimal morphing is 
considered as a numerical problem where a cost functional is minimized over 
a finite number of intermediate shapes. We introduce a theory of distortion 
minimal morphing over a continuous family of states in the context of morphs 
between $n$-dimensional oriented compact connected smooth manifolds without 
boundary embedded in $\R^{n+1}$ whose orientations are inherited from the 
usual orientation of $\R^{n+1}$. The natural cost functional (for stretching) measures the total relative 
change of volume with respect to a family of diffeomorphisms that defines the 
morph.  This functional is invariant under compositions with volume preserving 
diffeomorphisms; hence, the corresponding minimal morphs are not unique. On the other hand, we prove that the extremals of our functional are (strong) minima.  
Our main result is the existence of a distortion minimal morph 
(with respect to stretching) between every pair of isotopic submanifolds.

\section{Minimal distortion diffeomorphisms}
In this section we prove the existence of distortion minimal diffeomorphisms between diffeomorphic  $n$-dimensional oriented manifolds $M$ and $N$ (which are not necessarily embedded in $\R^{n+1}$)  with respective volume forms
$\w_M$ and $\w_N$.

Recall that the Jacobian of a diffeomorphism  $h: M\to N$ is defined by the equation
\[h^\ast \w_N = J(\w_M, \w_N)(h) \,\w_M,\]
where $h^\ast \w_N$ denotes the pullback of the volume form  $\w_N$ on $N$ by the diffeomorphism $h$ (see~\cite{AMR}). The Jacobian $J(h):=J(\w_M, \w_N)(h)$ depends on the diffeomorphism
and the volume forms.

The distortion (due to stretching) $\xi(m)$ at $m\in M$, with respect
 to a diffeomorphism $h:M \to N$, is defined by
\begin{equation}
\xi(m)=\lim \limits_{\ep \to 0} 
\frac{\big|\int_{h(A_\ep)} \w_N\big|-\big|\int_{A_\ep} \w_M \big| }
{\big| \int_{A_\ep} \w_M \big|}=\big|J(h)(m)\big|-1,
\end{equation}
where $A_\ep \subset M$,   for $\ep>0$,  is a nested family of (open) 
neighborhoods of the point $m \in M$ such that $A_\alpha\subseteq A_\beta$ 
whenever $\alpha>\beta>0$ and $\cap_{\ep>0} A_\ep=m$. 

In other words, 
the distortion is the infinitesimal relative change of volume with respect 
to $h$. It is easy to see that the definition of distortion does not
depend on the family of nested sets $A_\ep$.

We denote the set of all $C^2$ diffeomorphisms between
manifolds $M$ and $N$ by $\diff(M,N)$. 
The \emph{total distortion} functional $\Phi:\diff(M,N)\to \R$, 
with respect to the oriented 
manifolds $(M,\w_M)$ and
$(N,\w_N)$, is defined by
\begin{equation}\label{eq:func}
\Phi(h)=\int_M\Big(\big|J(h)(m)\big|-1\Big)^2\,  \w_M.
\end{equation}

We will establish necessary and sufficient conditions for a diffeomorphism 
$h:M \to N$ to be a minimum of the functional $\Phi$. Also, we will show that
a minimum always exists in $\diff(M,N)$ provided that the manifolds are
compact, connected, and without boundary.

As a useful notation, we let $\Chi(M)$ denote the set of smooth vector fields 
on the manifold $M$. Also, we recall a basic fact from global nonlinear analysis:
\emph{$\diff(M,N)$ is a Banach manifold and its tangent space at $h\in\diff(M,N)$
can be identified with $\Chi(N)$. }
Indeed, an element of $T_h \Diff(M,N)$ is an equivalence class of curves 
$[h_\ep]$, 
represented by a family of diffeomorphisms $h_\ep$ with $h_0=h$, where two curves
passing through $h$ are equivalent if they have the same derivative at $h$.
For each $n\in N$, this family defines a curve $\epsilon\mapsto h_\epsilon(h^{-1}(n))$ 
in $N$ that passes through $n$ at $\ep=0$; hence, 
it defines a vector $Y\in T_nN$ by 
\[Y(n):=\frac{d}{d\ep} h_\epsilon(h^{-1}(n))\Big |_{\ep=0}.\]
The vector field $Y\in \Chi(N)$ is thus associated with the equivalence 
class $[h_\ep]$. In fact, the vector field $Y$ does not depend on the choice of
the representative of the equivalence class.  
On the other hand, for $Y \in \Chi(N)$ with flow $\phi_t$, we associate
the curve $h_t=\phi_t \circ h$ in $\Diff(M,N)$. 
The (tangent) equivalence class of this
curve is an element in $T_h \Diff(M,N)$.

\begin{proposition}[Euler-Lagrange Equation]
\label{pr:SteadyMap}
Suppose that $M$ and $N$ are smooth  connected compact orientable manifolds 
without boundary. A $C^2$ diffeomorphism $h:M \to N$ is a critical point of
 the total distortion functional $\Phi$ if and only if $J(h)$ is constant.
\end{proposition}

\begin{proof}
Let $h_\ep:(-1,1) \to \diff(M,N)$ be a curve of diffeomorphisms from $M$ to $N$
such that $h_0=h$.  By definition,  $h \in \diff(M,N)$ is a critical point of the functional $\Phi(h)$, if 
$\frac{d}{d t}\Phi(h_t)\Big|_{t=0}=0$. Using the easily derived formula  \[\Phi(h)=\int_M J(h)^2 \w_M-2\vol(N)+\vol(M), \]
we note that $h$ is a critical point if and only if
$$
2\int_M J(h) \frac{d}{dt} J(h_t)|_{t=0}\, \w_M=0.
$$
Moreover, using the calculus of differential forms (see \cite{AMR} and note in particular that $L_Y$ is used to denote the Lie derivative in the direction of the vector field $Y$), 
we have that for $h^t=\psi_t \circ h$, where $\psi_t$ is the flow of $Y \in
\Chi(N)$,
\begin{eqnarray*}
\frac{d}{dt}(J(\psi_t\circ h) \w_M)\Big|_{t=0}
&=& \frac{d}{dt} \big((\psi_t\circ h)^*\w_N\big)\Big|_{t=0}\\
&=& h^* \frac{d}{dt} (\psi_t^*\w_N)\Big|_{t=0}\\
&=& h^*  \psi_t^*L_Y \w_N\big|_{t=0}\\
&=& h^*L_Y \w_N\\
&=& h^*(\div Y \w_N)\\
&=& (\div Y)\circ h J(h) \w_M.
\end{eqnarray*}

By Stokes' theorem and the properties of the $\wedge$-antiderivations $d$ and
$i_Y$, we have that
\begin{eqnarray*}
\frac{d}{d t}\Phi(h_t)\Big|_{t=0}&=&\int_M J(h)^2 \div Y \circ h \,\w_M
=\int_N J(h)\circ h^{-1}\div Y \w_N\\
&=&\int_N J(h)\circ h^{-1} L_Y \w_N
=\int_N J(h)\circ h^{-1} d\,i_Y \w_N\\
&=&
\int_N d(J(h)\circ h^{-1}\wedge i_Y \w_N)
-\int_N d(J(h)\circ h^{-1}) \wedge i_Y \w_N\\
&=&-\int_N i_Y (d(J(h)\circ h^{-1})\wedge \w_N)
+\int_N i_Y\big(d(J(h)\circ h^{-1})\big)\,\w_N\\
&=&
\int_N d(J(h)\circ h^{-1})(Y)\,\w_N.
\end{eqnarray*}
Hence, $h \in \diff(M,N)$ is a critical point of the functional $\Phi(h)$ if
and only if
$$
\int_N d(J(h)\circ h^{-1})(Y)\,\w_N=0
$$
for all $Y \in \Chi(N)$. It follows that if $J(h)$ is constant, then $h$ is a critical point of $\Phi$. 

To complete the proof  it suffices to show that if 
\begin{equation}
\label{eq155}
\int_N df(Y)\,\w_N=0
\end{equation}
for all $Y\in\Chi(N)$, then $df=0$, where $f:=(J(h)\circ h^{-1})$. 

Suppose, 
on the contrary, that there exists a continuous vector field $Y \in \Chi(N)$ such that $df(Y)(n) \neq 0$ for some point $n \in N$. Without
loss of generality,  we assume the inequality $df(Y)(n) > 0$. The map $df(Y):N\to \R$ is continuous. 
Therefore, there exists an open neighborhood $U \subset N$ of
the point $n \in N$ so that $df(Y)(p)>0$ for every $p \in U$. Using the standard bump function argument 
(see \cite{AMR}), we can construct a vector field $Z \in \Chi(N)$ supported in
$U$ such that 
$
\int_N df(Z) \w_N=\int_U df(Z) \w_N>0,
$
in contradiction to equality~\eqref{eq155}.
Hence, $df=0$. 
\end{proof}

\begin{definition}
A function $h \in \diff(M,N)$ is called a  \emph{\J map} if it is a critical point of
the total distortion functional $\Phi$. 
\end{definition}

As an immediate corollary of proposition~\ref{pr:SteadyMap}, we have
the following theorem.

\begin{theorem}
\label{JC}
A function $h \in \diff(M,N)$ is a \J map if and only if $J(h)$ is the 
constant function with value $\vol ( N)/\vol (M)$.
\end{theorem}

We will use the elementary properties of \J maps stated in the following lemma. The proof is left to the reader.

\begin{lemma}
\label{pr:compositions}
Compositions and inverses of \J maps  are \J maps.
\end{lemma}

Also we will use (the strong form) of Moser's theorem  on volume forms, which we state here for the convenience of the reader (see~\cite{M}).

\begin{theorem}
\label{MoserThm}
Let  $\tau_t$ be a family of volume forms defined for $t \in [0,1]$   on a compact manifold $M$. If 
\begin{equation}
\int_c \tau_t =\int_c \tau_0
\end{equation}
for every $n$-cycle $c$ on $M$,  then there exists a one-parameter family of 
diffeomorphisms $\phi_t:M\to M$ such that
\begin{equation}
\phi_t^\ast\tau_t=\tau_0
\end{equation}
and $\phi_0$ is the identity mapping. Moreover, the dependence of $\phi_t(m)$ 
on $m\in M$ and $t\in [0,1]$ is as smooth as in the family $\tau_t$.
\end{theorem}

\begin{theorem}\label{th:minphi}
If $(M,\w_M)$ and $(N,\w_N)$  are diffeomorphic 
$n$-dimensional compact 
connected oriented manifolds without boundary,  then
(i) there is a \J map from $M$ to $N$, (ii) every \J map from $M$ to $N$ 
minimizes the functional $\Phi$, and   
(iii) the minimum value of $\Phi$ is
\begin{equation}
\label{PhiMin}
\Phi_{min}=\frac{\big( \vol(M)-\vol(N)\big)^2}{\vol(M)}.
\end{equation}
\end{theorem}
\begin{proof}
To prove (i), choose a diffeomorphism $h \in \diff(M,N)$ and note that
the differential form $h^\ast \w_N$ is a volume on $M$. 
Define a new volume on $M$ as follows:
\[\bar{\w}_M=\frac{\vol(M)}{\int_M h^\ast \w_N} h^\ast \w_N.\]
Since \[\int_M \bar{\w}_M=\int_M {\w}_M\] and $M$ is compact, 
by an application of  Moser's theorem~\ref{MoserThm}, 
there exists a
 $C^2$ diffeomorphism $f:M \to M$ such that $\w_M=f^\ast \bar{\w}_M$. Hence,
\[\frac{\int_M h^\ast \w_N}{\vol(M)} \omega_M= 
(h \circ f)^\ast \w_N;\]
and, since $M$ is connected, we conclude that $J(h \circ
f)=\vol(N)/\vol(M)$ is constant and $DJ(h \circ f)=0$. 
Thus, $k=h \circ f$ is a \J map. 

To prove parts (ii) and (iii), note that 
if $k$ is an arbitrary \J map from $M$ to $N$, then
\[
\Phi(k)=\big(|J(k)|-1\big)^2 \vol(M)=
\frac{\big( \vol(M)-\vol(N)\big)^2}{\vol(M)}.\]
We claim that this value of $\Phi$ is its minimum.

Let $g \in \diff(M,N)$.
By the Cauchy-Schwartz inequality,
\begin{eqnarray*}
\Phi(g)&=&\int_M \big(|J(g)|-1\big)^2 \w_M\\ 
&\geq& \frac{1}{\vol(M)} \Big( \int_M \big(|J(g)|-1\big) \w_M \Big)^2\\
&=& \frac{1}{\vol(M)} \big( \vol(M)-\vol(N)\big)^2\\&=&\Phi(k),
\end{eqnarray*}
as required.
\end{proof}

\begin{example}
Let $S_r$ and $S_R$ be two-dimensional
round spheres of radii 
$r$ and $R$ (respectively)  centered at the origin in $\R^3$. 
Define $h:S_r \to S_R$ by
$h(p)=R/r \, p$ for $p=(x,y,z) \in S_r $. We will show that $h$ is a \J map.

Let $\omega_r $ (respectively, $w_R$)  be the standard volume forms on $S_r$ (respectively, $S_R$) generated
by the usual Riemannian metric $g_p(X,Y)=\langle X,Y \rangle$ for 
$X,Y \in \R^3$.

Using  the parametrization of $S_r$ and $S_R$ by spherical coordinates, it is
easy to compute that the Jacobian   $$J(\om_{r},
\om_{R})(h)(m)=R^2/r^2={\vol(S_R)}/{\vol(S_r)}$$ for all $m\in S_r$; 
hence, by theorem~\ref{th:minphi}, $h$ is a \J map. 
\end{example}

\begin{remark}[Harmonic maps]
For $h\in \Diff(M,N)$, the distortion functional~\eqref{eq:func} has value
\[
\Phi(h)=\int_M |J(h)|^2 \w_M-2 \vol(N)+\vol(M).
\]
Thus, it suffices to consider the minimization problem for the reduced functional
$\Psi$ given by
 \[
\Psi(h)=\int_M |J(h)|^2 \w_M.
\]
We note that if $M$ and $N$ are one-dimensional, then $\Psi$ is the same as
\[
\Psi(h)=\int_M |Dh|^2 \w_M.
\]
An extremal of this functional is called a harmonic map (see~\cite{EL, EL2, E}). Thus,
for the one-dimensional case, distortion minimal maps and harmonic maps coincide.
On the other hand, there seems to be no obvious relationship in the general case. 
\end{remark}

\section {Morphs of embedded manifolds}

We will discuss a minimization problem for morphs of compact connected 
boundaryless oriented n-dimensional smooth manifolds embedded in $\R^{n+1}$.
 
\subsection{Pairwise minimal morphs}
\begin{definition}
\label{df:morph}
Let $M$ and $N$ be compact connected oriented $n$-dimensional smooth 
manifolds without boundary embedded in
$\R^{n+1}$. A $C^1$ function
$H:[0,1] \times M \to \R^{n+1}$ is a morph from $M$ to $N$ if the
following conditions hold:
\begin{itemize}
\item[(i)] $p \mapsto H(t,p)$ is a diffeomorphism onto its image for each $t
\in I=[0,1]$;
\item[(ii)]  the image $M^t=H(t,M)$ is an $n$-dimensional manifold possessing all the
properties of $M$ and $N$ mentioned above;
\item[(iii)]  $p \mapsto H(0,p)$ is a diffeomorphism of $M$;
\item[(iv)]  the image of the map $p \mapsto H(1,p)$ is $N$.
\end{itemize}
We denote the set of all morphs between manifolds $M$ and $N$ by $\SM(M,N)$.
\end{definition}

For simplicity, we will consider only morphs $H$ such that $p \mapsto H(0,p)$ is the identity map. We assume that each manifold $M^t=H(t,M)$ (with  $M^0=M$ and $M^1=N$) is equipped with the volume form $\om_t=i_{\eta_t}\Omega$,
where $$\Omega=dx_1 \wedge dx_2 \wedge \ldots \wedge dx_{n+1}$$ is the standard
volume form on $\R^{n+1}$ and $\eta_t:M^t \to \R^{n+1}$ is the outer unit normal vector field on $M^t$ with respect to the usual metric on $\R^{n+1}$. Also,  as a convenient notation, we use $h^t=H(t,\cdot):M \to M^t$.
\begin{definition}
A morph $H$ is \emph{distortion pairwise minimal}
 (or, for brevity, \emph{pairwise minimal}) if $h^{s,t}=h^t \circ (h^s)^{-1}:M^s
\to M^t$ is a \J map for every $s$ and $t$. We denote the set of all distortion
pairwise minimal morphs between manifolds $M$ and $N$ by $\SP \SM(M,N)$. 
\end{definition}

By proposition \ref{pr:SteadyMap} and theorem \ref{th:minphi},
a morph $H$ is pairwise minimal if and only if each Jacobian
$J(\om_s, \om_t)(h^{s,t})$ is constant.

\begin{proposition}
\label{pr:pairwise}
Let $M=M^0$ and $N=M^1$ be $n$-dimensional manifolds as in definition~\ref{df:morph}
equipped with the (respective) volume forms $\om_0$ and $\om_1$. A morph $H$ between $M$ and $N$ is distortion pairwise minimal 
if and only if
\begin{equation}
\label{ConstPairwise}
\frac{J(\om_0, \om_t)(h^t)(m)}{\vol(M^t)}=\frac{1}{\vol(M)}
\end{equation}
for all $t \in [0,1]$ and $m\in M$.
\end{proposition}
\begin{proof}
Using lemma \ref{pr:compositions} and theorem 
\ref{th:minphi}, it suffices to prove that each map $h^t:M \to M^t$ is minimal 
if and only if the map~\eqref{ConstPairwise} is
constant. An application of theorem~\ref{JC} finishes the proof.
\end{proof}

\begin{proposition}
\label{pr:pwminimal}
Let $M$ and $N$ be $n$-dimensional manifolds as in proposition~\ref{pr:pairwise}. If there is a morph $G$ from $M$ to $N$, then there is a distortion pairwise
minimal morph between $M$ and $N$.
\end{proposition}

\begin{proof}
Fix a morph $G$ from $M$ to $N$ with the corresponding family of
diffeomorphisms $g^t:=G(t,\cdot)$, let $M^t:=G(t,M)$, 
and  consider the family of volume forms
$$
\bar{\om}_t=\frac{\vol(M)}{\vol(M^t)} (g^t)^\ast \om_t
$$
defined for $t\in [0,1]$.
It is easy to see that 
$$
\int_M \bar{\om}_t=\int_M \bar{\om_0};
$$
hence, by  Moser's theorem~\ref{MoserThm}, there is a family of diffeomorphisms 
$\alpha^t$ on $M$ such that $\alpha^t$ depends continuously on $t$ and
$\om_M=(\alpha^t)^\ast\bar{\om}_t$. It follows that 
$$(g^t \circ \alpha^t)^\ast
\om_t=\frac{\vol(M^t)}{\vol(M)} \om_M;$$ 
therefore, 
$$J(\om_M,\om_t)(g^t \circ
\alpha^t)(m)=\frac{\vol(M^t)}{\vol(M)}$$ 
for all $m\in M$. The morph $H$ corresponding
to the family $h^t:=g^t \circ
\alpha^t$ is the desired  distortion pairwise minimal morph.
\end{proof}

\subsection{Minimal morphs}
We will define distortion minimal morphs between embedded connected oriented
$n$-dimensional smooth manifolds without boundary.

For a morph $H$ from $M$ to $N$, let 
$E_{s,t}$ denote the total distortion of $h^{s,t}:M^s\to M^t$. We have that
\begin{eqnarray*}
E_{s,t}&=&\int_{M^s}\Big( \big|J(h^{s,t})\big|-1\Big)^2 \om_s\\
&=&\int_M\Big(\frac{J(h^t)}{J(h^s)}-1\Big)^2 J(h^s) \om_M.
\end{eqnarray*}
If $H$ is a $C^2$ morph, then $E_{s,t}$ is twice continuously differentiable with respect to $s$. By Taylor's theorem,  $E_{s,t}$ has the representation
$$
E_{s,t}=E_{t,t}+\frac{d}{ds}(E_{s,t})\big|_{s=t}(s-t)+
\frac{1}{2}\frac{d^2}{ds^2}(E_{s,t})\big|_{s=t}(s-t)^2+ O\big((s-t)^3\big).
$$
Note that $E_{t,t}$ and $\frac{d}{ds}(E_{s,t})\big|_{s=t}$
both vanish, and 
\begin{equation*}
\frac{1}{2}\frac{d^2}{ds^2}(E_{s,t})\big|_{s=t}
=\int_M \frac{\Big(\frac{d}{dt} J(h^t) \Big)^2}{J(h^t)} \om_M.
\end{equation*} 

\begin{definition}
The \emph{infinitesimal distortion} of a $C^2$ morph $H$ from $M$ to $N$ at 
$t\in [0,1]$ is 
$$
\ep^H(t)= \lim \limits_{s \to t} \frac{E_{s,t}}{(s-t)^2}=
\int_M \frac{\Big(\frac{d}{dt} J(h^t) \Big)^2}{J(h^t)} \om_M.
$$
The \emph{total distortion functional} $\Phi$ defined on such morphs is given by \begin{equation}
\label{eq4:1}
\Phi(H)=\int_0^1 \ep^H(t) dt= \int_0^1 
\Bigg(\int_M \frac{\Big(\frac{d}{dt} J(h^t) \Big)^2}{J(h^t)}
\om_M\Bigg)dt.
\end{equation}
\end{definition}

\begin{definition}
A $C^2$ morph is called a \emph{distortion minimal extremal} if it is an 
extremal of the functional
$\Phi$ with respect to $C^2$ morphs. 
A $C^2$ morph is called a \emph{distortion minimal morph} if it minimizes the 
functional $\Phi$.
\end{definition}

Note that $\ep^H(t)$ depends continuously on $t$ provided that $H$ is a $C^2$ 
morph.

\begin{lemma}
\label{Lemma:MinPM}
For every morph $H \in \SM(M,N)$ there exists a pairwise minimal morph 
$G \in \SP \SM(M,N)$ such that $\Phi(G)\leq\Phi(H)$.
In particular, if $H \in \SM(M,N)$ is a distortion minimal $C^2$ morph, then there exists a 
pairwise minimal morph $G \in \SP \SM(M,N)$ such that $\Phi(H)=\Phi(G)$.
\end{lemma}
\begin{proof}
Let $H \in \SM(M,N)$ be a morph with the intermediate states
$M^t=H(t,M)$. By proposition \ref{pr:pwminimal}, there exists a pairwise minimal
morph $G \in \SP\SM(M,N)$ with the same intermediate states. 
The deformation energy of transition maps satisfies the inequality
 $E_{s,t}(H)\geq E_{s,t}(G)$ for all $s,t \in [0,1]$ 
because $G$ is pairwise minimal.
Therefore, 
$\ep^H(t) \geq \ep^G(t)$ for all $t \in [0,1]$, and, consequently,
\begin{equation}
\label{eqw223}
\Phi(H) \geq \Phi(G)
\end{equation}
as required. 

If $H$ is distortion minimal,  the inequality $\Phi(H)\leq\Phi(G)$
holds.
Comparing the latter inequality with inequality
\eqref{eqw223}, we conclude that $\Phi(H)=\Phi(G)$.
\end{proof}

\begin{corollary}
\label{Cor:PairwInf}
\begin{itemize}
\item[(i)] The following inequality holds:
\begin{equation}
\inf_{G \in \SP\SM(M,N)} \Phi(G) \leq \inf_{H \in \SM(M,N)} \Phi(H).
\end{equation}
\item[(ii)]
If there exists a minimum $F$ of the total distortion functional $\Phi$ over 
the class $\SP\SM(M,N)$, then $F$ minimizes the functional $\Phi$ over the class 
$\SM(M,N)$ as well:
\begin{equation}
\Phi(F)=\min_{G \in \SP\SM(M,N)} \Phi(G) =\min_{H \in \SM(M,N)} \Phi(H).
\end{equation}
\end{itemize}
\end{corollary}

\begin{lemma} 
\label{Lemma:pairwEn}
The total distortion of a $C^2$ pairwise minimal morph
$H$ from $M$ to $N$ is
\begin{equation}
\label{eq4:2}
\Phi(H)=\int_0^1 \frac{\big(\frac{d}{dt}\vol(M^t)\big)^2}{\vol(M^t)}\, dt.
\end{equation}
\end{lemma}
\begin{proof}
The proof is an immediate consequence of formula \eqref{eq4:1} and 
proposition \ref{pr:pairwise}.
\end{proof}

\begin{lemma}
\label{Lemma:auxfunct}
Consider an auxiliary functional
\begin{equation}
\label{auxfct}
\Psi(\phi)=\int_0^1 \frac{\dot{\phi}^2}{\phi}\, dt
\end{equation}
defined on the admissible set $$Q=\big\{\phi \in C^1\big([0,1];\R_+\big):
\phi(0)=\vol(M), \phi(1)=\vol(N)\big\}.$$ 
The  functional~\eqref{auxfct} attains its minimum at $\phi \equiv \vol(M)$ if 
$\vol(M)=\vol(N)$ and at 
$\phi(t)=\Big[(\sqrt{\vol(M)} -\sqrt{\vol(N)}\,)t-\sqrt{\vol(M)}\,\Big]^2$ whenever
the volumes of $M$ and $N$ are distinct. The minimal value of the functional
$\Psi$ is
\begin{equation}
\min_{\rho \in Q} \Psi(\rho)=4 \big(\sqrt{\vol(N)}-\sqrt{\vol(M)}\big)^2.
\end{equation}

\end{lemma}
\begin{proof}The proof is a simple application of the Euler-Lagrange equation and 
the Cauchy-Schwarz inequality.
\end{proof}

Using corollary \ref{Cor:PairwInf} and lemma \ref{Lemma:auxfunct}, 
we will minimize the total distortion energy 
functional $\Phi$ over the set $\SM(M,N)$ of all morphs.

\begin{theorem}
Let $M$ and $N$ be two $n$-dimensional manifolds satisfying the assumptions 
of definition~\ref{df:morph}. If $M$ and $N$ are connected by a $C^2$ morph,  
 then there exists a distortion minimal morph. The minimal value of
$\Phi$ is
\begin{equation}
\label{PhiMin11}
\min_{H \in \SM(M,N)} \Phi(H) = 4
\big(\sqrt{\vol(N)}-\sqrt{\vol(M)}\big)^2.
\end{equation}

\end{theorem}
\begin{proof}
Let $G$ be a morph between $M$ and $N$. Without loss of generality, 
we assume that $G$ is pairwise minimal (see proposition \ref{pr:pwminimal}). Set
$$H(t,m)=\la(t)G(t,m),$$ where $\la:[0,1] \to \R$ is to be determined.

Note that if $M^t=H(t,M)$ and $W^t=G(t,M)$,  then
$$
\vol(M^t)=\int_M(h^t)^\ast\om_M=\big[\la(t)\big]^n \int_M(g^t)^\ast \om_M=
\big[\la(t)\big]^n \vol(W^t).
$$
Let 
$\phi(t)$ be the minimizer of the auxiliary functional $\Psi$ from lemma 
\ref{Lemma:auxfunct}, and 
define
$$\la(t)=\Big[\frac{\phi(t)}{\vol(W^t)}\Big]^{\frac{1}{n}}.$$ 
The corresponding volume
$\vol(M^t)=\phi(t)$; therefore, by corollary \ref{Cor:PairwInf} and 
lemma \ref{Lemma:auxfunct}, the morph 
$H$ minimizes the total distortion functional $\Phi$ over the class
$\SM(M,N)$ and $\Phi(H)=4 \big(\sqrt{\vol(N)}-\sqrt{\vol(M)}\big)^2$.
\end{proof}

The next result provides a basic class of distortion minimal morphs.

\begin{proposition}
Suppose that $M$ is an $n$-dimensional manifold embedded in $\R^{n+1}$ that satisfies the assumptions of definition~\ref{df:morph}. If $\al$ is a positive real number and  
 \[N:=\{\al m:  m \in M\},\] then 
\begin{itemize}
\item[(i)] $N$ is a manifold satisfying all of the assumptions of definition~\ref{df:morph}.  
\item[(ii)] the morph  given by the family of maps $h^t(m)=\la(t) m$, where
\[\la(t)=\vol(M)^{-\frac{1}{n}}
\Big[(\sqrt{\vol(M)} -\sqrt{\vol(N)}\,)t-\sqrt{\vol(M)}\,\Big]^{\frac{2}{n}},
\]
is distortion minimal.
\end{itemize}
\end{proposition}
\begin{proof}
Define $h^t(m)=\la(t) m$. It is easy to check that $h^t$ defines a morph from
$M$ to $N$. Also, we have that $J(h^t):=J(\om_M, \om_t)(h^t)=
\big[ \la(t) \big]^n$. 
Since $J(h^t)$ is constant on $M$, the family $h^t$ defines  
a pairwise minimal morph $H$.  

We will determine $\la(t)$  so that the morph $H$ becomes a minimizer of $\Phi$ over the class $\SM(M,N)$. Indeed, 
by lemma \ref{Lemma:auxfunct}, it suffices to choose $\la$ so that
$$
\vol(M^t)=\big[ \la(t) \big]^n\vol(M)=
\Big[(\sqrt{\vol(M)} -\sqrt{\vol(N)}\,)t-\sqrt{\vol(M)}\,\Big]^2,
$$
which  yields
$$
\la(t)=\vol(M)^{-\frac{1}{n}}
\Big[(\sqrt{\vol(M)} -\sqrt{\vol(N)}\,)t-\sqrt{\vol(M)}\,\Big]^{\frac{2}{n}}.
$$
The corresponding morph $H(t,m)=\la(t)m$ satisfies the equality
$$
\Phi(H)=\min_{G \in \SM(M,N)} \Phi(G).
$$
\end{proof}

\end{document}